\documentclass{article}
\usepackage{graphicx}
\usepackage{amssymb,amsmath,amsfonts}

\def\be{\begin{eqnarray}}
\def\ee{\end{eqnarray}}

\begin{document}

\thispagestyle{empty}

\hfill ITEP/TH-29/08

\bigskip

\centerline{\Large{Resultants and Contour Integrals}}

\bigskip

\centerline{\emph{A.Morozov and Sh.Shakirov}}

\bigskip

\centerline{ITEP, Moscow, Russia}

\centerline{MIPT, Dolgoprudny, Russia}

\bigskip

\centerline{ABSTRACT}

\bigskip

{\footnotesize
Resultants are important special functions used in description of non-linear phenomena. Resultant $R_{r_1, \ldots, r_n}$ defines a condition of solvability for a system of $n$ homogeneous polynomials of degrees $r_1, \ldots, r_n$ in $n$ variables, just in the same way as determinant does for a system of linear equations. Unfortunately, there is a lack of convenient formulas for resultants when the number of variables is large. In this paper we use Cauchy contour integrals to obtain a polynomial formula for resultants, which is expected to be useful in applications.
}

\bigskip

\bigskip

It is well-known, that a system of $n$ linear equations in $n$ variables has a non-zero solution if and only if the corresponding determinant vanishes. Solutions to non-homogeneous linear systems and Gaussian integrals (associated with linear equations of motion) are then expressed through determinants. Similarly in non-linear case a system of homogeneous polynomials $f_1( {\vec x} ), \ldots, f_n( {\vec x} )$ of degrees $r_1, \ldots, r_n$ in $n$ variables ${\vec x} = (x_1, \ldots, x_n)$ has a non-zero common root, if and only if coefficients of the system satisfy a constraint
\begin{align*}
R_{r_1, \ldots, r_n} \big\{ f_1, \ldots, f_n \big\} = 0
\end{align*}
where expression $R$ is called a \emph{resultant} of this system (the subscripts $r_1, \ldots, r_n$ will be often omitted in what follows). For linear systems resultant is just a determinant, the well-known and well-understood object of linear algebra. Unfortunately, the knowledge and understanding of general resultants is not equally deep, yet. Their study is a relatively new subject, which is natural to call \emph{non-linear algebra} \cite{GKZ, NOLINAL}.

 In analogue with linear algebra, solutions to non-homogeneous non-linear systems are expressed through resultants by generalized Craemer rule \cite{NOLINAL}, also resultants control singularities of non-Gaussian integrals. For this reason, resultants have all chances to become the central special functions of emerging non-linear physics. They are already used in physical applications -- see, for example, \cite{Appl1, Appl2, Appl3}.

 Practical application of non-linear algebra to scientific and engineering problems depends on development of efficient methods to handle the resultants -- as efficient as those for determinants in linear algebra. Like with determinants, not explicit expressions are needed (which take hundreds of pages already for low values of $n$ and $r$'s) but their properties and efficient computational algorithms, allowing to evaluate resultants in every particular situation. Several algorithms for calculation of resultants are already known \cite{Algo, Riem, LogDet}. In this paper we describe yet another algorithm, which can be well computerized and used in practice.

Our approach is based on the recursive relation
\begin{align}
\boxed{
\ \ \ \log R \big\{ f_1, f_2, \ldots, f_n \big\} = \oint \ldots \oint \log f_1 \ d\log f_2 \ \wedge \ldots \wedge \ d\log f_n + r_1 \log R \big\{ f_2,\ldots,f_n \big\} \Big|_{x_1 = 0} \ \ \
\label{M0}
}
\end{align}
where the contour integral is $(n-1)$-fold, functions $f_i$ are taken with arguments $ f_i(1,z_1,\ldots,z_{n-1}) $ and the contour of integration encircles infinity. A sketchy proof of this relation is given in Appendix A. The relation of resultants to contour integrals of this kind is also studied in \cite{Riem}, where promising connections with the theory of Riemann surfaces are pointed out. An explicit formula for $\log R$ follows after several iterations of (\ref{M0}): for example, for $n = 2$
$$ \log R_{r_1 r_2} \left\{ f_1, f_2 \right\} = \oint \log f_1(1,z) d \log f_2(1,z) + r_1 \log f_2(0,1)$$
for $n = 3$
\begin{align*}
\log R_{r_1 r_2 r_3} \left\{ f_1, f_2, f_3 \right\} =
& \oint \oint \log f_1(1,z_1,z_2) d \log f_2(1,z_1,z_2) \wedge d \log f_3(1,z_1,z_2) + \emph{}
\\ & r_1 \oint \log f_2(0,1,z) d \log f_3(0,1,z) + r_1 r_2 \log f_3(0,0,1)
\end{align*}
and so on. In this way, logarithm of the resultant is expressed through the simple contour integrals. However, we would like to calculate resultants, not just their logarithms. To obtain $R$ from $\log R$, we exponentiate the latter: for example, for $n = 2$
\be
 R_{r_1 r_2} \left\{ f_1, f_2 \right\} = \exp \left( \oint \log f_1(1,z) d \log f_2(1,z) + r_1 \log f_2(0,1) \right)
\label{M11}
\ee
Using this formula, it is straightforward to reproduce the well-known expression of the resultant of two polynomials through their roots: if
$$ f_1(1,z) = a \prod\limits_{i = 1}^{r_1} (z - \alpha_i) $$

$$ f_2(1,z) = b \prod\limits_{i = 1}^{r_2} (z - \beta_i) $$
then
$$ R_{r_1 r_2} \left\{ f_1, f_2 \right\} = \exp \left( \sum\limits_{i = 1}^{r_1} \sum\limits_{j = 1}^{r_2} \oint \log (z - \alpha_i) d \log (z - \beta_j) + r_2 \log a + r_1 \log b \right) = a^{r_{2}} b^{r_{1}} \prod\limits_{i = 1}^{r_1} \prod\limits_{j = 1}^{r_2} ( \beta_{j} - \alpha_i ) $$
Exponential formula (\ref{M11}) and its direct analogues for $n > 2$ have several advantages. The most important are simplicity and explicitness: such formulas have a transparent structure, which is easy to understand and memorize. Another one is universality: they can be written for any $n$, not just for low values. However, there is one serious drawback as well. Resultant is known to be a homogeneous polynomial of degree $$ d_i = \dfrac{r_1 r_2 \ldots r_n}{r_i} $$ in coefficients of the $i$-th equation $f_i$, and of total degree $$ d = d_1 + \ldots + d_n = r_1 r_2 \ldots r_n \left( \dfrac{1}{r_1} + \dfrac{1}{r_2} + \ldots + \dfrac{1}{r_n} \right) $$ in coefficients of all equations, but in exponential formula the polynomial nature of the resultant is not explicit. Therefore, it is not fully convenient for practical calculations of resultants (although it is conceptually adequate and can be used in theoretical considerations).

\smallskip

We would like to convert the exponential formula into another form, which would be explicitly polynomial. Following \cite{LogDet}, we suggest to apply (\ref{M0}) to a \emph{shifted} system
\[ \left\{ \begin{array}{c}
{\tilde f}_1( {\vec x}) = (x_1)^{r_1} - \lambda_1 f_1( {\vec x}) \\
\noalign{\medskip}{\tilde f}_2( {\vec x}) = (x_2)^{r_2} - \lambda_2 f_2( {\vec x}) \\
\noalign{\medskip}\ldots\\
\noalign{\medskip}{\tilde f}_n( {\vec x}) = (x_n)^{r_n} - \lambda_n f_n( {\vec x}) \\
\end{array} \right. \]
and expand the shifted logarithms in the integrands into Taylor series in powers of the "spectral parameters" $\lambda_i$. We obtain the following series expansion
\be
\log R \big\{ {\tilde f}_1, \ldots, {\tilde f}_n \big\} = - \sum\limits_{k_1 = 0}^{\infty} \ldots \sum\limits_{k_n = 0}^{\infty} T_{k_1 k_2 \ldots k_n}\big\{ f_1, \ldots, f_n \big\} \cdot \lambda_1^{k_1} \lambda_2^{k_2} \ldots \lambda_n^{k_n}
\label{MTr}
\ee
(a minus sign is just a convention) where particular Taylor components $T_{k_1 k_2 \ldots k_n}$ are homogeneous polynomial expressions of degree $k_i$ in coefficients of $f_i$. We call them \emph{traces} \cite{LogDet} of a non-linear system $f_1( {\vec x} ), \ldots, f_n( {\vec x} )$. Explicit formulas for traces are obtained in Appendix B, by making a series expansion in (\ref{M0}) and carrying out the remaining contour integrals.

By exponentiating, we obtain the shifted resultant:
$$ R \big\{ {\tilde f}_1, \ldots, {\tilde f}_n \big\} = \exp \left( - \sum\limits_{k_1 = 0}^{\infty} \ldots \sum\limits_{k_n = 0}^{\infty} T_{k_1 k_2 \ldots k_n} \cdot \lambda_1^{k_1} \lambda_2^{k_2} \ldots \lambda_n^{k_n} \right)$$
The original resultant $ R \big\{ f_1, \ldots, f_n \big\} $ is equal to $(-1)^{d}$ times the coefficient of $ \lambda_1^{d_1} \lambda_2^{d_2} \ldots \lambda_n^{d_n} $ in $R \big\{ {\tilde f}_1, \ldots, {\tilde f}_n \big\} $. We can extract it, expanding the right hand side in powers of $\lambda_i$:
$$  \exp \left( - \sum\limits_{k_1 = 0}^{\infty} \ldots \sum\limits_{k_n = 0}^{\infty} T_{k_1 k_2 \ldots k_n} \cdot \lambda_1^{k_1} \lambda_2^{k_2} \ldots \lambda_n^{k_n} \right) = \sum\limits_{k_1 = 0}^{\infty} \ldots \sum\limits_{k_n = 0}^{\infty} {\cal P}_{k_1 k_2 \ldots k_n} \cdot \lambda_1^{k_1} \lambda_2^{k_2} \ldots \lambda_n^{k_n}$$
Such power series expansion of $\exp \left( S(x) \right)$, where $S(x)$ is itself a power series (perhaps, of many variables), is often called a \emph{Schur expansion}. It implies the following relation between Taylor components:
$$
{\cal P}_{k_1 k_2 \ldots k_n} \ = \ \sum\limits_{m = 1}^{k_1 + \ldots + k_n} \ \dfrac{(-1)^{m}}{m!} \ \sum\limits_{{\vec v}_1 + {\vec v}_2 + \ldots + {\vec v}_m = {\vec k}} \ T_{{\vec v}_1} \ T_{{\vec v}_2} \ \ldots \ T_{{\vec v}_m}
$$
where the sum is taken over all ordered partitions of a a vector ${\vec k} = (k_1, k_2, \ldots, k_n)$ into $m$ vectors, denoted as ${\vec v_1}, \ldots, {\vec v_m}$. An ordered partition is a way of writing a vector with integer components as a sum of vectors with integer components, where the order of the items is significant. Polynomials ${\cal P}$ are often called multi-Schur polynomials. For example,
\[
\begin{array}{cc}
{\cal P}_{1, 0}  = - T_{1, 0}\\
\\
{\cal P}_{2, 0}  = - T_{2, 0} + T_{1, 0}^2/2\\
\\
{\cal P}_{2, 1}  = - T_{2,1} + T_{2,0} T_{0,1} + T_{1,0} T_{1,1} - T_{1,0}^2 T_{0,1}/2\\
\\
{\cal P}_{2, 1, 0} = - T_{2,1,0} + T_{2,0,0} T_{0,1,0} + T_{1,0,0} T_{1,1,0} - T_{1,0,0}^2 T_{0,1,0}/2\\
\\
{\cal P}_{1, 1, 1} = - T_{1,1,1} + T_{1,0,0} T_{0,1,1} + T_{0,1,0} T_{1,0,1} + T_{1,1,0} T_{0,0,1} - T_{0,1,0} T_{0,0,1} T_{1,0,0}\\
\end{array}
\]
Resultant is a certain multi-Schur polynomial of traces:
\begin{align}
\boxed{
\ \ R \big\{ f_1, \ldots, f_n \big\} = (-1)^d \ {\cal P}_{d_1,\ldots,d_n} = \sum\limits_{m = 1}^{d} \ \dfrac{(-1)^{m + d}}{m!} \ \sum\limits_{{\vec v}_1 + {\vec v}_2 + \ldots + {\vec v}_m = {\vec d}} \ T_{{\vec v}_1} \ T_{{\vec v}_2} \ \ldots \ T_{{\vec v}_m} \ \
\label{M1}
}
\end{align}
The explicit formula for traces, obtained in Appendix B, has a form
\begin{align}
\boxed{
\ \ T_{\vec k} = T_{k_1 k_2 \ldots k_n} = \dfrac{1}{k_1 k_2 \ldots k_n} \cdot \sum\limits_{ r_{ij} } \ \det\limits_{2 \leq i,j \leq n} \big( \delta_{ij} r_{i} k_{i} - r_{ij} \big) \prod\limits_{i = 1}^{n} (f_i)^{k_i}_{r_{i1}, r_{i2}, \ldots ,r_{in}} \ \
}
\label{M2}
\end{align}
where
\begin{align*}
 (f)^{k}_{j_1,j_2,\ldots, j_n} = \mbox{ coefficient of } x_1^{j_1} x_2^{j_2} \ldots x_n^{j_n} \mbox{ in } f(x_1,x_2,\ldots,x_n)^k
\end{align*}
and the sum goes over all non-negative integer $n \times n$ matrices $r_{ij}$ such that the sum of entries in any $i$-th row or $i$-th coloumn is fixed and equals $r_{i} k_{i}$. Formula (\ref{M2}) is valid for positive $k_1, \ldots, k_n$. The cases when some $k_i = 0$ are equally simple, but they are somewhat degenerate and considered separately in the Appendix B.

Together, formulas (\ref{M1}) and (\ref{M2}) constitute our final result, an explicit polynomial formula for the resultant of $n$ homogeneous polynomials $f_1, \ldots, f_n$ in $n$ variables. They provide a rather fast and effective computer algorithm for explicit evaluation of resultant for a finite set of polynomials $f_i$. A MAPLE program, based on this algorithm, is presented in the appendix C.

\section*{Appendix A: Proof of the relation (\ref{M0}) }

We begin from the well-known Poisson product formula \cite{GKZ, LogDet}. It states, that the resultant of a system
\be
\left\{
\begin{array}{lll}
f_1 \left( x_1, x_2, \ldots, x_n \right) = 0\\
\noalign{\medskip}f_2 \left( x_1, x_2, \ldots, x_n \right) = 0\\
\noalign{\medskip}\ldots \\
\noalign{\medskip}f_n \left( x_1, x_2, \ldots, x_n \right) = 0\\
\end{array}
\right.
\label{M4}
\ee
is equal to the product
\begin{align}
R \{ f_1, f_2, \ldots, f_n \} = C \cdot f_1 \left( {\vec \Lambda^{(1)}} \right) f_1 \left( {\vec \Lambda^{(2)}} \right) \ldots f_1 \left( {\vec \Lambda^{(N)}} \right)
\label{M5}
\end{align}
where vectors ${\vec \Lambda^{(i)}}$ are all the $N = r_2 r_3 \ldots r_n$ common roots of polynomials $f_2, \ldots, f_n$:
$$
\left\{
\begin{array}{lll}
f_2 \left( {\vec \Lambda^{(i)}} \right) = 0\\
\noalign{\medskip}\ldots \\
\noalign{\medskip}f_n \left( {\vec \Lambda^{(i)}} \right) = 0\\
\end{array}
\right.
$$
This formula holds for the following reason. By definition of solvability, (\ref{M4}) is solvable iff $f_1$ vanishes on (at least) one vector ${\vec \Lambda^{(i)}} \neq 0$. What is the same, resultant of (\ref{M4}) vanishes iff a product
$$f_1 \left( {\vec \Lambda^{(1)}} \right) f_1 \left( {\vec \Lambda^{(2)}} \right) \ldots f_1 \left( {\vec \Lambda^{(N)}} \right)$$ vanishes. Since resultant is a polynomial, it must be divisible by this product. The constant of proportionality $C$ depends on the normalization of roots (but not on $f_1$) and we assume that ${\vec \Lambda^{(i)}}$ are normalized in such a way, that $C = 1$. We proceed by rewriting this formula as
$$ R \{ f_1, f_2, \ldots, f_n \} = \left( \Lambda^{(1)}_{1} \Lambda^{(2)}_{1} \ldots \Lambda^{(N)}_{1} \right)^{r_1} \cdot f_1 \left( {\vec \mu^{(1)}} \right) f_1 \left( {\vec \mu^{(2)}} \right) \ldots f_1 \left( {\vec \mu^{(N)}} \right)$$
where $ {\vec \mu^{(i)}} = {\vec \Lambda^{(i)}} / \Lambda^{(i)}_{1} $ are the re-normalized roots, which have the first component $\mu^{(i)}_1 = 1$. Since at the same time it follows from (\ref{M5}) that
$$ R \{ x_1, f_2, \ldots, f_n \} = \Lambda^{(1)}_{1} \Lambda^{(2)}_{1} \ldots \Lambda^{(N)}_{1} $$
we have
$$ R \{ f_1, f_2, \ldots, f_n \} = R \{ x_1, f_2, \ldots, f_n \}^{r_1} \cdot f_1 \left( {\vec \mu^{(1)}} \right) f_1 \left( {\vec \mu^{(2)}} \right) \ldots f_1 \left( {\vec \mu^{(N)}} \right)$$
Passing to logarithms, we obtain
\be
 \log R \{ f_1, f_2, \ldots, f_n \} = r_1 \log R \{ x_1, f_2, \ldots, f_n \} + \sum\limits_{ i = 1 }^{N} \log f_1 \left( {\vec \mu^{(i)}} \right)
\label{M7}
\ee
In this form, Poisson product formula admits a contour integral representation. This is because the second term is a sum of values of one function at roots of the other(s). Using Cauchy residue theorem, such sums can be conveniently represented as contour integrals. Consider first an example, familiar from the elementary theory of functions of one complex variable. Let $z$ be a complex variable, let $f(z)$ and $g(z)$ be polynomials of degrees $|f|$ and $|g|$. Then, the sum of values of $\log f$ at roots of $g$ is given by
$$ \sum\limits_{ i = 1 }^{|g|} \log f\left(\mu^{(i)}\right) = \oint \log f(z) d\log g(z) $$
where the contour encircles all the $|g|$ roots of $g$, denoted by $\mu^{(i)}, \ i = 1,2,\ldots,|g|$. This is obvious, because
$$ g(z) = g_{0} \prod\limits_{i = 1}^{|g|} \left(z - \mu^{(i)}\right) $$
and
$$ \oint \log f(z) d\log g(z) = \sum\limits_{i = 1}^{|g|} \oint \dfrac{\log f(z)}{z - \mu^{(i)}} \ dz = \sum\limits_{ i = 1 }^{|g|} \log f\left(\mu^{(i)}\right)$$
This example is so simple, because any function of one complex variable is factorizable into linear factors. Such factorization has no analogue in the case of several complex variables, therefore, the proof is slightly more complicated in multidimensional situation.

\paragraph{Proposition.} The second term in the right hand side of (\ref{M7}) is given by a $(n-1)$-fold contour integral
$$ \sum\limits_{ i = 1}^{N} \log f_1 \left( {\vec \mu^{(i)}} \right) = \oint \ldots \oint \log f_1(1,z_2,\ldots,z_n) \prod\limits_{k = 2}^{n} d\log f_k(1,z_2,\ldots,z_n) $$
where the contour encircles all the $N = r_2 r_3 \ldots r_n$ common roots of $f_2, \ldots, f_n$, denoted by ${\vec \mu}^{(i)}$.

\paragraph{Proof.} Making a non-linear change of variables
\[
\left\{
\begin{array}{ccc}
y_2 = f_2(1,z_2,\ldots,z_n)\\
\noalign{\medskip}y_3 = f_3(1,z_2,\ldots,z_n)\\
\noalign{\medskip}\ldots\\
\noalign{\medskip}y_n = f_n(1,z_2,\ldots,z_n)\\
\end{array}
\right.
\]
we obtain
$$ \oint \ldots \oint \log f_1(1,z_2,\ldots,z_n) \prod\limits_{k = 2}^{n} d\log f_k(1,z_2,\ldots,z_n) = \sum\limits_{ i = 1}^{N} \oint \ldots \oint \log f_1 \left(1,{\vec z}^{(i)}\left( {\vec y} \right) \right) d\log y_2 \ldots d\log y_n$$
Since $f_k$ are polynomials of degrees $r_k$, the inverse change of variables is multi-valued and has $N = r_2 r_3 \ldots r_n$ branches, denoted as
$ {\vec z}^{(i)}\left( {\vec y} \right), \ i = 1,2,\ldots,N $ in the above formula. The integral turns into a sum over branches. Integration over $y$-variables is trivial
$$ \oint \ldots \oint \log f_1 \left(1,{\vec z}^{(i)}\left( {\vec y} \right) \right) \dfrac{dy_2 \ldots dy_n}{y_2 \ldots y_n} = \log f_1 \left(1,{\vec z}^{(i)}\left( {\vec 0} \right) \right) $$
due to the Cauchy residue theorem. Obviously, the inverse change of variables maps (0,0,\ldots,0) into ${\vec \mu^{(i)}}$:
$$\oint \ldots \oint \log f_1(1,z_2,\ldots,z_n) \prod\limits_{k = 2}^{n} d\log f_k(1,z_2,\ldots,z_n) = \sum\limits_{ i = 1 }^{N} \log f_1 \left(1,{\vec z}^{(i)}\left( {\vec 0} \right) \right) = \sum\limits_{ i = 1}^{N} \log f_1 \left( {\vec \mu^{(i)}} \right) $$
The proposition is proved. Therefore, we have
\be
 \log R \{ f_1, f_2, \ldots, f_n \} = r_1 \log R \{ x_1, f_2, \ldots, f_n \} + \oint \ldots \oint \log f_1(1,z_2,\ldots,z_n) \prod\limits_{k = 2}^{n} d\log f_k(1,z_2,\ldots,z_n)
\label{M8}
\ee
Relation (\ref{M0}) follows from (\ref{M8}) and the simple identity
$$ R \big\{ x_1, f_2(x_1,x_2,\ldots,x_n), \ldots, f_n(x_1,x_2,\ldots,x_n) \big\} = R \big\{ f_2(0,x_2,\ldots,x_n), \ldots, f_n(0,x_2,\ldots,x_n) \big\} $$
implied by the definition of resultant as a solvability condition.

\section*{Appendix B: Explicit formula for the traces (\ref{MTr})}

We are going to show, that for positive $k_i$, the traces $T_{k_1 k_2 \ldots k_n}$ are given by
\[
\begin{array}{cccc}
\\ \nonumber
T_{k_1 k_2 \ldots k_n} = \sum\limits_{ r_{ij} } \ \det\limits_{2 \leq i,j \leq n} \big( \delta_{ij} r_{i} k_{i} - r_{ij} \big) \cdot \prod\limits_{i = 1}^{n} \dfrac{(f_i)^{k_i}_{r_{i1}, r_{i2}, \ldots ,r_{in}}}{k_i}
\\ \nonumber
\end{array}
\]
where
\begin{align*} (f)^{k}_{j_1,j_2,\ldots, j_n} = \mbox{ coefficient of } x_1^{j_1} x_2^{j_2} \ldots x_n^{j_n} \mbox{ in } f(x_1,x_2,\ldots,x_n)^k \end{align*} for any polynomial $f$, and the sum goes over all non-negative integer $n \times n$ matrices $r_{ij}$ such, that the sum of entries in any $i$-th row or $i$-th coloumn is fixed and equals $r_{i} k_{i}$: $$\sum\limits_{j = 1}^{n} r_{ij} = \sum\limits_{j = 1}^{n} r_{ji} = r_{i} k_{i}$$
Indeed, making  in (\ref{M0}) a series expansion
$$ \log \left( z_{i}^{r_i} - \lambda_i f_{i} \right) = r_i \log z_{i} - \sum\limits_{k = 1}^{\infty} \dfrac{1}{k} \ f_i^{k} z_{i}^{- r_i k} \lambda_i^{k} $$
for positive $k_i$ we have
\\
$$ T_{k_1,k_2,\ldots, k_n}(f) = \dfrac{ (-1)^{n+1} } {k_1 k_2 \ldots k_n } \oint \ldots \oint f_1^{k_1} d \left( f_2^{k_2} z_{2}^{- r_2 k_2} \right) \wedge d \left( f_3^{k_3} z_{3}^{- r_3 k_3} \right) \wedge \ldots \wedge d \left( f_n^{k_n} z_{n}^{- r_n k_n} \right) $$
\\
Using
$$
f_i^{k_i}(1,z_2,\ldots,z_n) = \sum\limits_{i_1 + \ldots + i_n = r_i k_i} (f_i)^{k_i}_{i_1 i_2 \ldots i_n} z_{2}^{i_2} \ldots z_{n}^{i_n}\\
$$
we obtain after taking derivatives
\\
$$ T_{k_1,k_2,\ldots, k_n}(f) = \sum\limits_{r_{ij}} \det\limits_{2 \leq i,j \leq n} \big( r_{ij} - \delta_{ij} r_{i} k_{i} \big) \prod\limits_{i = 1}^{n} \dfrac{ (f_i)^{k_i}_{r_{i1}, r_{i2}, \ldots ,r_{in}} } {k_i} \cdot \prod\limits_{i = 2}^{n} \oint z_{i}^{r_{i1} + \ldots + r_{in} - r_{i}k_{i} - 1} dz_{i}  $$
\\
The Cauchy residue theorem implies a selection rule $r_{i1} + \ldots + r_{in} = r_{i}k_{i}$. The transverse selection rule $r_{1i} + \ldots + r_{ni} = r_{i}k_{i}$ is implied by homogenity of polynomials $f_i$. Consider several examples in low dimensions:
\[
\begin{array}{ccc}
T_{k_1}(f) = \sum\limits_{r_{ij}} \dfrac{(f_1)^{k_1}_{r_{11}}}{k_1} = \dfrac{(f_1)^{k_1}_{r_1 k_1}}{k_1} \\
\\
T_{k_1, k_2}(f) = \sum\limits_{r_{ij}} \left( r_{2} k_{2} - r_{22} \right) \dfrac{(f_1)^{k_1}_{r_{11} r_{12}}}{k_1} \dfrac{(f_2)^{k_2}_{r_{21} r_{22}}}{k_2} \\
\\
T_{k_1, k_2, k_3} (f) = \sum\limits_{r_{ij}} \left|
\begin{array}{ccc}
r_{2} k_{2} - r_{22} & - r_{23} \\
\noalign{\medskip}- r_{32} & r_{3} k_{3} - r_{33}\\
\end{array}
\right| \cdot \dfrac{ (f_1)^{k_1}_{r_{11} r_{12} r_{13}} } {k_1 }  \dfrac{(f_2)^{k_2}_{r_{21} r_{22} r_{23}}} {k_2 } \dfrac{(f_3)^{k_3}_{r_{31} r_{32} r_{33}}} {k_3 }\\
\\
\\
T_{k_1, k_2, k_3, k_4} (f) = \\
\\
 \sum\limits_{r_{ij}} \left|
\begin{array}{ccc}
r_{2} k_{2} - r_{22} & - r_{23} & - r_{24} \\
\noalign{\medskip}- r_{32} & r_{3} k_{3} - r_{33} & - r_{34} \\
\noalign{\medskip}- r_{42} & - r_{43} & r_{4} k_{4} - r_{44} \\
\end{array}
\right| \cdot \dfrac{ (f_1)^{k_1}_{r_{11} r_{12} r_{13} r_{14}} } {k_1 } \dfrac{ (f_2)^{k_2}_{r_{21} r_{22} r_{23} r_{24}} } {k_2 } \dfrac{ (f_3)^{k_3}_{r_{31} r_{32} r_{33} r_{34}} } {k_3 } \dfrac{ (f_4)^{k_4}_{r_{41} r_{42} r_{43} r_{44}} } {k_4 }
\\
\end{array}
\]
\smallskip
The cases, when some $k_i$ are equal to zero, also make no difficulty. If some $k_i = 0$, then
$$T_{k_1,k_2,\ldots,k_n}\{ f_1,f_2,\ldots,f_n \} = r_i \left. T_{k_1,k_2,\ldots,k_{i-1},k_{i+1},\ldots,k_n}\{  f_1, f_2,\ldots, f_{i-1},f_{i+1},\ldots,f_{n}\} \right|_{x_i = 0}$$
where the trace in the right hand side is taken in variables $x_1, \ldots, x_{i-1}, x_{i+1}, \ldots, x_{n}$. It is an easy exercise to prove this statement, making a series expansion in (\ref{M0}).

\section*{Appendix C: MAPLE Program}

\begin{verbatim}

#################################### RESULTANT OF N HOMOGENEOUS POLYNOMIALS IN N VARIABLES
interface(warnlevel = 0):
interface(displayprecision = 2):
with(combinat):
with(linalg):

##################### INITIAL DATA

N := 3:

f[1] := a * x[1]^2 + alpha * x[2]*x[3];
f[2] := b * x[2]^2 + alpha * x[1]*x[3];
f[3] := c * x[3]^2 + alpha * x[1]*x[2];

###################################

BeginTime := time():

R := [seq(0, p = 1..N)]:

for i from 1 to N do
  if not ( select(t -> op(0,t) = x, indets(f[i])) subset {seq(x[p], p = 1..N)} ) then
    print( " Error: wrong variables. Must be ", seq(x[p], p = 1..N) );
    return Error;
  fi:
  R[i] := simplify( ln( subs( {seq(x[p] = exp(1) * x[p], p = 1..N)}, f[i] ) / f[i] ) ):
  if nops(indets(R[i])) <> 0 then
    print( " Error: non-homogeneous polynomials " );
    return Error;
  fi:
end do:

deg := [seq( mul (R[pp], pp = 1..N) / R[p], p = 1..N)]:
xDeg := add( deg[p], p = 1..N):

Decompose := (u,v) -> {seq( [seq( composition(u+v,v)[p][q] - 1,
q = 1..nops(composition(u+v,v)[p]))], p = 1..nops(composition(u+v,v)))}:

xDet := u -> if rowdim(u) = 0 then return 1: else return det(u): fi:

xTrace := proc(K):
  x1 := select( t -> K[t] = 0, [$1..N]):
  x2 := select( t -> K[t] <> 0, [$1..N]):
  for p from 1 to nops(x2) do
    F[p] := - subs({seq(x[ss] = 0, ss in x1), x[x2[1]] = 1},
    (f[x2[p]])^K[x2[p]]/x[x2[p]]^(R[x2[p]]*K[x2[p]])/K[x2[p]] ):
  end do:
  dat := mul(R[ss], ss in x1) * F[1] * xDet(matrix(nops(x2)-1,nops(x2)-1,
  (i,j) -> diff( F[i+1], x[x2[j+1]] ) ) );
  for p from 2 to nops(x2) do
    dat := eval(coeff( dat, x[x2[p]], -1));
  end do:
  return simplify(dat):
end proc:

xxTrace := kk -> add( xTrace(vv), vv in Decompose(kk, N) ):

xxSchur := kk -> add( mul( Tr[v[lll]], lll = 1..nops(v) )/nops(v)! *
multinomial(nops(v), seq( numboccur(v, convert(v,set)[i]), i = 1..nops(convert(v, set)))),
 v in partition(kk) ):

print( " Calculating resultant of type R[", op(R), "]" );
print( " Need to calculate ", xDeg, "traces " );
print( " ---------------------------------- " );

for kkk from 1 to xDeg do
  Tr[kkk] := xxTrace(kkk);
  print( kkk );
end do:

print( " ---------------------------------- " );

print( " Resultant = ", factor( simplify( (-1)^xDeg * xxSchur(xDeg) )));

print( " Done in ", time() - BeginTime, " sec. ");

\end{verbatim}
\section*{Acknowledgements}

It is a pleasure to thank V.Dolotin and all the participants of ITEP seminars on non-linear algebra for useful discussions. This work is partly supported by Russian Federal Atomic Energy Agency and Russian Academy of Sciences, by the joint grant 06-01-92059-CE, by NWO project 047.011.2004.026, by INTAS grant 05-1000008-7865, by ANR-05-BLAN-0029-01 project, by RFBR grant 07-02-00645 and by the Russian President's Grant of Support for the Scientific Schools NSh-3035.2008.2. The work of Sh.Shakirov is also partly supported by the Dynasty Foundation.

\end{document}